\documentclass[12pt]{amsart}
\usepackage{amsmath,amssymb,amsfonts,amscd,verbatim,graphicx}
\usepackage{xypic}
\newtheorem{theorem}[subsection]{Theorem}
\newtheorem{lemma}[subsection]{Lemma}

\newtheorem{claim}[subsection]{Claim}
\theoremstyle{definition}
\newtheorem{Definition}[subsection]{Definition}

\theoremstyle{remark}

\numberwithin{equation}{section}

\theoremstyle{definition}

\begin{document}

\title{A polynomial invariant and duality for triangulations}
\author{Vyacheslav Krushkal and David Renardy}
\thanks{This research was supported in part by the NSF}
\address{Department of Mathematics, University of Virginia, Charlottesville, VA 22904}
\email{krushkal\char 64 virginia.edu, dr7jx\char 64 virginia.edu}

\begin{abstract}
The Tutte polynomial ${T}_G(X,Y)$ of a graph $G$ is a classical invariant, important in combinatorics and
statistical mechanics. An essential feature of the Tutte polynomial is the duality for planar graphs $G$,
$T_G(X,Y)\; =\;  {T}_{G^*}(Y,X)$ where $G^*$ denotes the dual graph. We examine this property from the perspective
of manifold topology, formulating polynomial invariants for higher-dimensional simplicial complexes.
Polynomial duality for triangulations of a sphere follows as a consequence of Alexander duality.

The main goal of this paper is to introduce and begin the study of a more general $4$-variable polynomial for triangulations and handle decompositions of orientable manifolds.
Polynomial duality in this case is a consequence of Poincar\'{e} duality on manifolds.
In dimension $2$ these invariants specialize to the well-known polynomial invariants of ribbon graphs
defined by B. Bollob\'{a}s and O. Riordan.
Examples and specific evaluations of the polynomials are discussed.
\end{abstract}

\maketitle

\section{Introduction}
The Tutte polynomial is a classical invariant of graphs and matroids, important in combinatorics,
knot theory, and statistical mechanics \cite{Tutte1, Bollobas}.
Two properties of the Tutte polynomial of graphs are of particular interest: the contraction-deletion rule, and the duality
$T_G(X,Y)\; =\;  T_{G^*}(Y,X)$
where $G$ is a planar graph and $G^*$ is its dual. Several generalizations of the Tutte polynomial have been introduced for graphs embedded in surfaces. These polynomials reflect both the combinatorial properties of graphs and the topological invariants of their embeddings. This line of research was initiated by the work of B. Bollob\'{a}s and O. Riordan \cite{BR} with their introduction of a polynomial invariant of ribbon graphs. A further contribution was made by the first named author in \cite{Krushkal} where a four-variable generalization of the Tutte polynomial for graphs on surfaces was defined, satisfying a rather natural duality relation.

In recent years various combinatorial, topological and analytic invariants and properties of graphs have been extended to the setting of simplicial and CW complexes. The main purpose of this paper is to investigate and generalize the Tutte polynomial and its duality for ribbon graphs in the context of higher-dimensional complexes and their embeddings in manifolds.

We consider two versions of the Tutte polynomial for simplicial complexes, and more generally for
CW complexes. The first version, $T_K(X,Y)$,
is defined for an arbitrary simplicial (and more generally CW) complex $K$. If $K$ is a triangulation of the sphere $S^{2n}$, a direct
analogue of the Tutte duality, $T_K(X,Y)\; =\;  T_{K^*}(Y,X)$, where $K^*$ is the dual cell complex,
holds as a consequence of Alexander duality.
Recent results of \cite{BBC}, in particular showing that the evaluation
$T_K(0,0)$ gives the number of simplicial spanning trees in $K$, are mentioned in section \ref{spanning trees}.

The polynomial $T_K(X,Y)$ may be interpreted as the Tutte polynomial of a {\em matroid} associated to the simplicial chain
complex of $K$. The topological duality discussed above precisely corresponds to the notion of matroid duality.

A more general four-variable polynomial invariant is defined for a CW complex $K$ (of dimension $\geq n$) embedded in an orientable $2n$-dimensional manifold $M$, using the intersection pairing structure on the middle-dimensional homology group $H_n(M)$. If $K$ is
a triangulation of $M$, then as a consequence of the Poincar\'{e} duality on $M$ we prove
$$P_{K,M}(X,Y,A,B)=P_{K^*,M}(Y,X,B,A)$$
where $K^*$ is the dual CW complex. The polynomial $P$ is analogously defined for handle decompositions of $M$, and the duality stated above also holds
for a handle decomposition and its dual. While the polynomial $T_K$ discussed above may be defined in terms of a simplicial matroid, an interpretation of the invariant $P_{K, M}$ in the context of matroid theory is not currently known (except for the case of graphs on surfaces, see discussion in section \ref{questions}). If $M$ is oriented and its dimension is divisible by $4$, the polynomial $P_{K, M}$ may be further refined using the decomposition into positive-definite and negative-definite subspaces associated to the intersection pairing on $M$.

The polynomial $P$ for graphs on surfaces (corresponding to $n=1$) was introduced in \cite{Krushkal}. Its definition was motivated in part by 
questions in statistical mechanics, specifically the Potts model on surfaces, and applications to topology (the Jones polynomial for virtual knots). Moreover, it unified many previously defined
(and seemingly unrelated) invariants of graphs on surfaces: \cite{Krushkal}, \cite{BBC}, \cite{ACEMS} respectively established that the Bollob\'{a}s-Riordan polynomial \cite{BR}, the Bott polynomial \cite{Bott} and the Las Vergnas polynomial \cite{Las Vergnas} are in fact all specializations of the polynomial $P$. 

The purpose of this paper is to investigate the duality and other properties of $P$ for higher-dimensional complexes. One may view the polynomial $P$ (for graphs and in higher dimensions)
as an invariant of thickenings of a given complex $K$, taking into account the topological information about the thickening. $2$-dimensional thickenings of graphs are known as ribbon graphs, and more generally in higher dimensions a rigorous way to define a ``thickening'' of a complex $K$ is to consider a handle decomposition whose spine is $K$.
We give examples of the polynomial for specific handle decompositions of several manifolds. Since the topology and the combinatorics of complexes (and especially of their thickenings) of dimension $2$ and higher are substantially more involved than those in dimension $1$, and since many properties and applications of polynomial invariants of ribbon graphs have been discovered by various authors in recent years, we also expect a rich theory associated with the polynomial $P$ in higher dimensions.

The organization of the paper is as follows. After reviewing the background material in section \ref{background},
we formulate the polynomial $T$ and prove its duality for triangulations of a sphere in section \ref{sphere section}.
An interpretation of the polynomial $T$ in terms of simplicial matroids, and the relation between matroid duality and topological duality
are presented in section \ref{matroid section}.
The four-variable polynomial $P$ is defined, and the corresponding duality theorem is proved in section \ref{manifold polynomial}. Lemma \ref{relation with Tutte} establishes that the polynomial $T$ is
a specialization of the more general invariant $P$, where the topological information reflecting the embedding of the complex $K$ into a manifold $M$ is disregarded.
Section \ref{evaluation} contains calculations of the polynomial $P$ for specific manifolds. It also gives examples
of evaluations of the polynomials $T, P$ generalizing the classical fact that the number of spanning trees of a graph $G$ is the evaluation
of the Tutte polynomial $T_G(0,0)$ . Section \ref{generalizations} discusses generalizations of the polynomials
$T$ and $P$. The final section \ref{questions} mentions several questions motivated by our results.

\subsection*{Acknowledgements}
We are grateful to Sergei Chmutov for his comments on the paper and especially for bringing to our attention the recent results of \cite{BBC}.
We would like to thank Victor Reiner and Hugh Thomas for their comments on the original version of this paper.

We also would like to thank the referee for making a number of comments and suggestions which improved the exposition of the paper.

\section{Background} \label{background}
This section reviews some basic notions in topology and combinatorics that will be used throughout the paper. We begin by discussing the Tutte polynomial for graphs.

\subsection{Graphs and the Tutte Polynomial} \label{graph section}
A graph $G$ is defined by a collection of vertices, $V$, together with a specified collection $E$ of pairs of vertices called edges. A graph is {\em planar} if it is embedded into the plane. Given a planar graph $G=(V,E)$, one can construct its dual graph $G^*= (V^*,E^*)$,
whose vertices $V^*$ correspond to the connected components of ${\mathbb R}^2\smallsetminus G$. Two vertices in $V^*$ are connected by an edge in $E^*$ if the regions they represent are adjacent along an edge from $E$.

Given a graph $G=(V,E)$, a {\em spanning subgraph} $H\subset G$, $H=(V, E')$, has the same vertex set $V$ as $G$ and $E'\subseteq E$. Consider the following
normalization of the Tutte polynomial \cite{Tutte1}.

\begin{Definition} The Tutte Polynomial of a graph $G$ is defined by:
\begin{equation} \label{Tutte definition}
T_G(X,Y)=\sum_{H \subset G}{(X-1)^{c(H)-c(G)}(Y-1)^{n(H)}}
\end{equation}
where the summation is taken over all spanning subgraphs $H\subset G$ of $H$, $c(H)$ is the number of connected components of $H$, and $n(H)$ is the $nullity$ of $H$ given by the rank of the first homology group of $H$.
(The nullity may also be defined combinatorially as $n(H)=c(H)+|E(H)|-|V(H)|$).
We will also use a different normalization of this polynomial, known as the {\em rank-generating polynomial} (cf. \cite{Bollobas, BO}):
\begin{equation} \label{Rank definition}
R_G(X,Y)=\sum_{H \subset G}{X^{c(H)-c(G)}Y^{n(H)}}
\end{equation}
\end{Definition}
The Tutte polynomial  is a classical invariant in graph theory, and in particular its one variable specializations, the chromatic polynomial and the flow
polynomial, are of considerable independent interest.
While the Tutte polynomial encodes many properties of a graph, the main focus of this paper is on its duality relation, a feature important for applications to topology and statistical mechanics. Specifically, for a connected planar graph $G$ and its dual $G^*$, $T_G(X,Y)=T_{G^*}(Y,X)$ \cite{Bollobas}.

\subsection{Triangulations and Duality} \label{triangulation section}

Our goal is to define a version of the Tutte polynomial for simplicial complexes (and more generally, for cell complexes) and to establish a duality relation similar to that of planar graphs.
For this purpose, we consider triangulations and handle decompositions of smooth manifolds.
Planar graphs may be thought of as cellulations of the sphere $S^2$, and when thickened a planar graph provides a handle
decomposition of the sphere (vertices give rise to disks, edges to ribbons glued on along the boundaries of disks).
Handle decompositions for manifolds of higher dimensions provide a geometirc form of duality that generalizes dualization for planar graphs.
In the remainder of this section we will summarize the relevant material about triangulations and handle decompositions, a more detailed
account may be found in \cite{RS}, \cite{Milnor}.

\begin{Definition} A {\em triangulation} of a topological space $X$ is a simplicial complex $K$ homeomorphic to $X$ along with a homeomorphism $h:K\rightarrow X$.
\end{Definition}
In the formulation of the Tutte polynomial for graphs (\ref{Tutte definition}), the sum is taken over all spanning subgraphs. Framing this in the language of simplicial complexes, the sum is over all subcomplexes of top dimension $1$ such that the entire $0$-skeleton is included in each subcomplex. We generalize this condition to higher dimensions.

\begin{Definition}\label{spanning definition}
 Given a simplicial (or CW) complex $K$ of dimension $ \geq n$, let $L$ be its $n$-dimensional subcomplex.
Call $L$ a {\em spanning $n$-subcomplex} of $K$ if their $(n-1)$-skeletons coincide, $L^{(n-1)}=K^{(n-1)}$.
\end{Definition}

If $K$ is a finite complex (this will be the case throughout the paper), there are $2^{C_n}$ spanning $n$-subcomplexes of $K$, where $C_n$ is the number of $n$-cells of $K$.
Suppose $K$ is a triangulation of an $n$-manifold $M$. The dual cell complex $K^*$ (defined below) of $K$ has an important property that the $k$-cells of $K^*$ are in one-to-one correspondence with the $(n-k)$-cells of $K$. This correspondence has deeper consequences that will help establish the duality statements for the polynomials introduced in this paper. We will next discuss a combinatorial construction of $K^*$, equivalently the dual complex may also be
defined using dual handle decompositions, see section \ref{handle section}.

Let $M$ be an $n$-manifold and let $K$ be a triangulation of $M$. Given a $k$-simplex $\sigma _k$ of $K$, its dual, $D \sigma _k$, is an $(n-k)$-cell formed by taking the
union of all simplices of the barycentric subdivision  that contain the centroid of $\sigma _k$ as a vertex and that are transverse to $\sigma _k$.
Taking the collection of all such $D \sigma _k$ for $k=0,1,...,n$ gives the desired dual cell complex $K^*$.

For example, the dual of an $n$-simplex is the $0$-cell corresponding to the centroid. For a $0$-simplex, the dual is the union of all simplices in the barycentric subdivision of M that have that $0$-simplex as a vertex, which gives the $n$-cell that contains the $0$-cell as a centroid. Thinking of a planar graph as a cell decomposition of the sphere $S^2$, the dual graph corresponds exactly to the dual complex. For more on triangulations and the construction of the dual complex consult \cite{RS}. The geometric duality on manifolds is best understood
in the context of handle decompositions, discussed in the following subsection.
\subsection{Handle Decompositions} \label{handle section}

A handle decomposition of a manifold is analogous to a cell decomposition of a topological space. The goal is to understand the entire space as a union of $n$-balls pieced together by prescribed attaching maps. Let $0\leq k\leq n$ and consider the $n$-ball, $D^n$, as the product $D^k \times D^{n-k}$. It is attached to a manifold along a part of its boundary. More precisely:

\begin{Definition} Let $M$ be an $n$-manifold with boundary. Let $0\leq k \leq n$ and let $f\! : (\partial D^k)\times D^{n-k}\longrightarrow
\partial M$ be an embedding. Then $M\cup_f (D^k\times D^{n-k})$ is called the result of {\em attaching a $k$-handle to $M$}.
Some standard terminology:  $f$ is the {\em attaching map} of the handle, $f(\partial D^k \times  0 )$ is the {\em attaching sphere},
$f(D^k \times 0 )$ is the {\em core} and $f(0  \times D^{n-k} )$ is the {\em co-core} of the handle.
\end{Definition}

It is a basic and central fact in Morse theory \cite{Milnor} that any smooth manifold $M$ admits a handle
decomposition. In fact, given a triangulation $K$ of $M$, the simplices of $K$ may be thickened to produce
a handle decomposition (see \cite[p. 82]{RS}).
Conversely, a handle decomposition may be retracted to give a cell decomposition of $M$, see \cite[p. 83]{RS}.
(Each handle is retracted onto its core, using the product structure of the handle.)

Given a handle decomposition of $M$, each $k$-handle $D^k\times D^{n-k}$ dually may be thought of as
an $(n-k)$-handle attached along the complementary part of its boundary, $D^k\times \partial D^{n-k}$.
This gives rise to a dual handle decomposition of $M$. Given a triangulation $K$ of $M$, thickening $K$ gives rise
to a handle decomposition ${\mathcal H}$, dualizing gives a handle decomposition ${\mathcal H}^*$, then retracting
the handles of ${\mathcal H}^*$ onto their cores gives a complex $K^*$. This is a construction of the dual
complex $K^*$, alternative to the combinatorial construction discussed above.

\section{The Tutte polynomial for complexes and duality for triangulations of a sphere} \label{sphere section}

A natural generalization of the Tutte polynomial to higher dimensions defined below is formulated using
homology groups of subcomplexes of a given complex $K$. All homology groups considered in this paper are
taken with real coefficients, $H_i(\, .\, ; {\mathbb R})$. We refer the reader to \cite{Hatcher} as a basic reference in
homology theory. Denote by $|H_n(L)|$ the rank (dimension) of the $n$th homology group of $L$. The following
definition is formulated for CW complexes, but the reader interested in the more restricted class
of simplicial complexes may replace the term ``CW'' by ``simplicial'' and all definitions and proofs hold in this context as well.

\begin{Definition} \label{sphere definition}
Let $K$ be a CW complex of dimension $\geq n$. Define
\begin{equation} \label{general Tutte}
 T_K(X,Y)=\sum_{L \subset K^{(n)}} {X^{|H_{n-1}(L)|-|H_{n-1}(K)|}Y^{|H_n(L)|}}
\end{equation}
where the summation is taken over all  spanning $n$-subcomplexes $L$ of $K$ (see definition
\ref{spanning definition}).
\end{Definition}

A more precise notation for the polynomial defined in (\ref{general Tutte}) is $T_{K,n}$, including a reference to the dimension $n$. However in the case of main interest in this paper $K$ will be a triangulation of a $2n$-dimensional manifold and $n$ in definition (\ref{general Tutte}) will always be half the dimension of the ambient
manifold. Therefore the reference to $n$ is omitted from our notation.

This definition lends itself to a number of generalizations, for example see section \ref{other T}.
Section \ref{manifold polynomial} below defines a $4$-variable polynomial
$P_K$ for a complex $K$ embedded in a $2n$-manifold $M$, giving $T_K$ as a particular specialization.
If $K$ is a $1$-complex (i.e. a graph), the definition of the polynomial $T_K$ coincides with the classical rank-generating polynomial (\ref{Rank definition}), a renormalization of the Tutte polynomial (\ref{Tutte definition}). 
The contraction-deletion rule for $T_K$ is analyzed in \cite{BBC}.
The following is a
generalization of the duality for the Tutte polynomial of planar graphs.

\begin{theorem} \label{sphere duality}
{\sl Let $K$ be a triangulation of $S^{2n}$, then
\begin{equation} \label{duality equation}
T_{K^{(n)}}(X,Y)=T_{K^{*(n)}}(Y,X)
\end{equation} where $K^{(n)}$ is the $n$-skeleton of
$K$ and $K^{*(n)}$ is the $n$-skeleton of the dual complex $K^*$.}
\end{theorem}

The duality relation (\ref{duality equation}) also holds in a more general setting where $K$ is the CW complex (not necessarily
a triangulation) associated to a handle decomposition of $S^{2n}$, see section \ref{handle section}. The proof of the theorem still holds when $S^{2n}$ is
replaced by an orientable $2n$-manifold $M$ such that $H_{n-1}(M)=H_n(M)=0$. (For manifolds $M$ without
this vanishing condition on homology, the more general polynomial $P$ of section \ref{manifold polynomial} provides the right context for the duality statement.)

A generalization of the polynomial $T$, taking into account the cardinality of the torsion subgroups of the homology (with $\mathbb Z$ coefficients) of the subcomplexes $L$, has been suggested
in \cite{BBC}. It is shown in \cite{BBC} that this refinement still satisfies the duality analogous to (\ref{duality equation}), see section \ref{spanning trees}  for further
discussion of this invariant.

The definition (\ref{general Tutte}) of $T_K$ is a special case of the rank-generating/Tutte polynomial of a matroid, see section \ref{matroid section}, associated to the simplicial/cellular matroid of $K$. Theorem \ref{sphere duality} follows from matroid duality for the chain complexes associated to a triangulation and its dual. This fact has been considered in the literature, cf. \cite[Proposition 6.1]{DKM1},  see section \ref{matroid section} for further details. We give a  more topological argument below as a warm-up for the proof of a more general duality statement for the polynomial $P$ defined in section \ref{manifold polynomial}.

\begin{proof}[Proof of theorem \ref{sphere duality}]
When $K$ is a triangulation of the sphere, $H_{n-1}(K)=H_{n-1}(S^{2n})=0$, therefore definition \ref{sphere definition} in this context reads
\begin{equation} \label{general Tutte for triangulation}
 T_K(X,Y)=\sum_{L \subset K^{(n)}} {X^{|H_{n-1}(L)|}Y^{|H_n(L)|}}
\end{equation}
Given a spanning $n$-subcomplex $L\subset K$, let $\overline L$ be the spanning $n$-subcomplex of the dual complex $K^*$ containing all of the $n$-cells of $K^*$ except those dual to the $n$-simplices of $K$. An important ingredient of the proof is the observation that $\overline  L$ is homotopy equivalent to $S^{2n}\smallsetminus L$. More generally:

\begin{lemma}\label{useful lemma} {\sl Let $M$ be a closed orientable $2n$-manifold, and let $K$ be a triangulation of $M$. Let $L$ be a spanning $n$-subcomplex of $K$ and let $\overline L$ be the corresponding $n$-subcomplex of $K^*$ described above. Then $\overline L$ is homotopy equivalent to $M\smallsetminus L$.}
\end{lemma}
\begin{proof}
Recall from \cite{RS} and section \ref{handle section} above  that triangulations give rise to handle decompositions. Specifically, construct a handle decomposition $\mathcal{H}$ of $M$  by thickening each $k$-simplex in $K$ to a $k$-handle, $k=0,1,\ldots, 2n$. Consider all handles that result from thickening $L \subset K$ and call this collection $\mathcal{H}_L$.  Notice that the handles in $\mathcal{H}\smallsetminus \mathcal{H}_L$ are thickenings of the simplices in $K\smallsetminus L$, and these are precisely the simplices of $K$ dual to those in $\overline L$.
Considering these handles dually, $\mathcal{H}\smallsetminus \mathcal{H}_L$ is a thickening of $\overline L$ (and hence it is homotopy equivalent to $\overline L$). To summarize,
$M={\mathcal H}_L\cup ({\mathcal H}\smallsetminus {\mathcal H}_L)$, where ${\mathcal H}_L$ is homotopy equivalent to $L$  (in fact $L$ is a deformation retract of ${\mathcal H}_L$) and ${\mathcal H}\smallsetminus {\mathcal H}_L$ is homotopy equivalent to $\overline L$. It follows that $M\smallsetminus L$ is homotopy equivalent to $M\smallsetminus \mathcal{H}_L$ and this in turn is homotopy equivalent to $\overline L$, finishing the proof of lemma \ref{useful lemma}.
\end{proof}

For each spanning $n$-subcomplex $L\subset K$, consider the corresponding $\overline L\subset K^*$  as above.  Observe that
$$|H_{n-1}(L)|=|H_{n}(\overline L)|.$$
Indeed, one has $|H_n(X)| = |H^n(X)|$ for any topological space $X$, and Alexander duality \cite{Hatcher} for the sphere states that
\[H^{n-1}(L)\cong H_{2n-n}(S^{2n}\smallsetminus L)=H_{n}(S^{2n}\smallsetminus L).\]
Since $\overline L$ is homotopy equivalent to $S^{2n}\smallsetminus L$ we conclude that $|H_{n-1}(L)|=|H_n(\overline L)|$. By the symmetry of our construction, this also
gives $|H_{n-1}(\overline L)|=|H_n(L)|$.
Since the spanning subcomplexes $L, \overline L$ are in $1-1$ correspondence,  the corresponding terms in the expansion (\ref{general Tutte for triangulation}) of the two sides of
(\ref{duality equation}) are equal. This concludes the proof of theorem \ref{sphere duality}.
\end{proof}

\section{The simplicial Matroid} \label{matroid section}

A matroid is a finite set with a notion of independence that generalizes the concept of linear independence in vector spaces. This notion was introduced by H. Whitney \cite{Whitney}, detailed expositions may be found in \cite{Oxley}, \cite{Tutte4}, \cite{Welsh}.
We begin with some background definitions for matroids and the formulation of the Tutte polynomial in this context.
\begin{Definition}
A {\em matroid} is a finite set $E$ with a specified collection $I$ of subsets of $E$, called the {\em independent sets} of $E$, such that:

(1) \parbox[t]{14.25cm}{$\emptyset\in I$.}\\
(2) \parbox[t]{14.25cm}{If $B\in I$ and $A\subset B$ then $A\in I$.}\\
(3) \parbox[t]{14.25cm}{If $A,B \in I$ and $|A|>|B|$ then there exists $a\in A\smallsetminus B$ such that
$a\cup B \in I$.}

A maximal independent set in $E$ is called a {\em basis} for the matroid.
\end{Definition}
An important example is given by {\em graph matroids}. A finite graph $G$ gives rise to a matroid as follows: take the set of all edges to be the set $E$ and call a collection of edges independent if and only if it does not contain a cycle. Equivalently, the matroid associated to a graph may be defined using the (adjacency) linear map from
the vector space spanned by its edges to the one spanned by its vertices. (This is a basic example of a simplicial matroid, and this point of view is examined in more detail further below.)
\begin{Definition}
If a set $E$ with independent sets $I$ forms a matroid, then a {\em rank function} $r$ assigns a non-negative integer to every subset of $E$ such that:

(1) \parbox[t]{14.25cm}{ $r(A) \leq |A|$ for all $A\subset E$. (Here $|A|$ denotes the cardinality of $A$.)}\\
(2) \parbox[t]{14.25cm}{If $A\subset B \subset E$, then $r(A) \leq r(B).$}\\
(3) \parbox[t]{14.25cm}{If $A,B \subset E$ then $r(A \cup B) + r(A\cap B) \leq r(A) + r(B)$.}
\end{Definition}

A matroid is determined by its rank function. One could alternatively define the independent sets of $E$ as the sets $A\subset E$ with $|A|=r(A)$.

\subsection{Duality on Matroids}

There is a natural notion of duality for matroids. If $M=(E,I)$ is a finite matroid, then the dual matroid $M^*$ is obtained by taking the same underlying set $E$ and the condition that a set, $A$, is a basis in $M^*$  if and only if $E\smallsetminus A$ is a basis in $M$. An important result of Kuratowski gives as a corollary that for a graphic matroid $M$, the dual matroid is graphic if and only if $M$ is the matroid of a planar graph. The rank function of the dual matroid is given by $r^*(A)= |A|-r(E) +r(E\smallsetminus A)$.

The Tutte polynomial of a matroid $M=(E,I)$ with rank function $r$ is defined as follows:
\[ T_M(X,Y)=\sum_{A \subset E}{(X-1)^{r(E)-r(A)}(Y-1)^{|A|-r(A)}}\]
As in section \ref{graph section}, we also consider a renormalization, the rank-generating polynomial \cite{BO}
$R_M(X,Y)=\sum_{A \subset E}{X^{r(E)-r(A)}Y^{|A|-r(A)}}.$
\begin{theorem} \label{dual matroid Tutte} \cite{Welsh}
{\sl The Tutte polynomial for matroids satisfies the duality $$T_M(X,Y)=T_{M^*}(Y,X).$$}
\end{theorem}

\subsection{An alternative proof of Theorem \ref{sphere duality}}
An alternative argument relies on a construction of a matroid whose rank-generating  polynomial coincides with the polynomial defined in (\ref{general Tutte}). The {\em simplicial} matroid has been considered by a number of authors, and the duality statement follows for example from \cite[Proposition 6.1]{DKM1}. We include this discussion since it may be less familiar to the topologist reader, and also it may be of interest since a matroid interpretation of the polynomial $P$ in section \ref{manifold polynomial} is not currently known, see item 4 in section \ref{questions}.

Let $K$ be as in the statement of theorem \ref{sphere duality}, and consider the simplicial chain complex for $K$:
\[\begin{CD}... @>>> C_{n+1}(K) @>{\partial_{n+1}}>> C_n(K) @>{\partial_n}>> C_{n-1}(K) @>>> ...
\end{CD}\]
Recall that $C_i(K)$ is the free abelian group generated by the $i$-simplices of $K$.
\begin{Definition} \label{simplicial matroid}
Given a simplicial complex $K$ of dimension $\geq n$, let $E$ be the set of $n$-simplices of $K$ (also thought of as a specific choice of generators of $C_n(K)$). A collection of elements of $E$ is said to be independent if their images under $\partial_n $ are linearly independent in $C_{n-1}(K)$. The resulting matroid $M(K)=(E,I)$ is called the {\em simplicial matroid} associated to $K$.
(For a CW complex and the corresponding cellular chain complex, this matroid will
be referred to as the {\em cellular matroid} associated to $K$.)
\end{Definition}
As in definition \ref{sphere definition}, a more precise notation for this matroid is $M(K,n)$, however $n$ should be clear from the context and is omitted from the notation.
The rank function $r$  on the matroid $M(K)$ is defined by $r(A)= {\rm rank}\, (\partial_n(A))$ where $\partial_n(A)$ is the subgroup of $C_{n-1}(K)$ generated by the images of the elements of $A$ under ${\partial}_n $.

This matroid has been studied by a number of authors, see \cite{Cordovil}, \cite{CrapoRota}. We will investigate the simplicial matroid in the context of triangulations of a sphere, showing that matroid duality then precisely corresponds to topological duality.

\begin{lemma}\label{dual lemma}  {\sl
(1) Let $K$ be a simplicial complex of dimension $\geq n$. Then the polynomial $T_K$ defined
in (\ref{general Tutte}) coincides with the rank-generating polynomial $R_{M(K)}$  associated to the simplicial matroid $M(K)$.

(2) If $K$ is a triangulation of the sphere $S^{2n}$ then the dual matroid $(M(K))^*$ coincides with the cellular matroid associated to the dual cell complex $K^*$, that is $(M(K))^*=M(K^*)$. }
\end{lemma}

\begin{proof}[Proof of lemma \ref{dual lemma}.] Note that subsets $A\subset E$ correspond to spanning $n$-subcomplexes of $K$ (definition \ref{spanning definition}).
Given $A\subset E$, consider the corresponding $n$-subcomplex $L$ (equal to the $(n-1)$-skeleton of $K$ union with the $n$-simplices corresponding to the elements of $A$). We get the following commutative diagram induced by the inclusion $L\subset K$:
\[\begin{CD} 0  @>{\partial_{n+1}^L}>> C_n(L) @>{\partial_n^L}>> C_{n-1}(L) @>{\partial_{n-1}^L}>> C_{n-2}(L) @>>>...\\
@VVV @VVV @| @| \\
C_{n+1}(K) @>{\partial_{n+1}^K}>> C_n(K) @>{\partial_n^K}>> C_{n-1}(K) @>{\partial_{n-1}^K}>>  C_{n-2}(K) @>>>...
\end{CD}\]
where $C_{n+1}(L)= 0$ since $L$ has no $(n+1)$-cells. Thus $H_n(L)\cong {\rm ker} \, \partial_n^L$, and $|H_n(L)| = |{\rm ker}\, \partial_n^L| = |A|-r(A)$.

Note that $|H_{n-1}(L)|-|H_{n-1}(K)|=|{\rm ker}\, \partial_{n-1}^L/{\rm im}\, \partial_n^L| - |{\rm ker}\, \partial_{n-1}^K/{\rm im}\, \partial_n^K| = |{\rm im}\, \partial_n^K|- |{\rm im}\, \partial_n^L| = r(E)-r(A)$.
Thus the rank-generating  polynomial $R_{M(K)}$ of the matroid $M(K)=(E,I)$ coincides with the polynomial $T_K$ defined in (\ref{general Tutte}).

We will now show that when $K$ is a triangulation of $S^{2n}$, the dual matroid $(M(K))^*$ coincides with the simplicial matroid structure described above applied to $K^*$. Consider the simplicial chain complex for $K^*$:
\[\begin{CD}... @>>> C_{n+1}(K^*) @>{\partial^*_{n+1}}>> C_n(K^*) @>{\partial_n^*}>> C_{n-1}(K^*) @>>> ...
\end{CD}\]
From the construction of the dual cell complex $K^*$ (sections \ref{triangulation section}, \ref{handle section}) it is clear that $C_i(K) \cong C_{2n-i}(K^*)$, and moreover $\partial^*_{i}$ is the adjoint of $\partial_{2n-i+1}$.

By definition $A$ is a basis of $M(K)$ if and only if $E\smallsetminus A$ is a basis for $(M(K))^*$.
It follows from the properties of the adjoint map that $\partial _n^K(A)$ is a linear basis for ${\rm im}\, \partial_n^K$ if and only if $\partial_n^*(E^*\smallsetminus A^*)$ is a basis for
${\rm im}\, \partial_n^*$.
Thus $(M(K))^*=(E^*,I^*)$ where a set $A^*\subset E^*$ is in $I^*$ iff $\partial^*_n(A^*)$ is linearly independent in $C_{n-1}(K^*)$. Therefore $(M(K))^*= M(K^*)$, concluding the proof of lemma \ref{dual lemma}. 
\end{proof}

It follows from the lemma that if $K$ is a triangulation of $S^{2n}$ then $T_{(M(K))^*}=T_{M(K^*)}$. This observation together with theorem \ref{dual matroid Tutte} gives an alternative proof of theorem \ref{sphere duality}.

\section{A Polynomial invariant for Triangulations of an Orientable Manifold} \label{manifold polynomial}

Let M be a closed oriented $2n$-dimensional manifold. Let $K$ be a simplicial (or CW) complex embedded in $M$, for example a triangulation of $M$, and let $L$ be a spanning $n$-subcomplex of $K$ (see definition \ref{spanning definition}) with $i\! : L \longrightarrow M$ being the embedding. Recall that throughout this paper all homology groups are taken with coefficients in ${\mathbb R}$. Define:
\begin{equation} \label{kernel}
k(L)={\rm rank}\, ({\rm ker}\,(i_*: H_n(L) \rightarrow H_n(M)))
\end{equation}

Let $\cdot$ denote the intersection pairing on $M$: $$ \mathbf{\cdot} \, : H_n(M) \times H_n(M) \rightarrow \mathbb{R}.$$ 
(The intersection pairing is the Poincar\'{e} dual of the cup product in cohomology \cite[p. 249]{Hatcher}.)
Consider the following vector spaces defined using the intersection pairing:
\begin{equation} \label{image}
V=V(L)={\rm image}\, (i_*: H_n(L) \rightarrow H_n(M)),
\end{equation}
\begin{equation} \label{perp}
V^{\perp}=V^{\perp}(L)= \{ u\in H_n(M) | \, \forall v\in V(L), u\cdot v=0 \}
\end{equation}
Consider two invariants of the embedding $L\longrightarrow M$:
\begin{equation} \label{s}
s(L):= {\rm dim}\, (V/(V\cap V^{\perp}))\; \; \, {\rm  and} \; \; \, s^{\perp}(L):= {\rm dim}\, (V^{\perp}/(V\cap V^{\perp})).
\end{equation}
This construction is motivated by the work in \cite{Krushkal} corresponding to the case $n=1$.
In the case $n=1$ ($M$ is a surface and $L$ is a graph in $M$) there is a geometric interpretation of the invariants $s$, $s^{\perp}$: $s$ equals twice the genus of the surface obtained as the regular neighborhood of the graph $L$ in $M$, and similarly $s^{\perp}$ is twice the genus of the regular neighborhood of the dual graph.

Another invariant of the embedding $L \longrightarrow M$ is
\begin{equation} \label{lagrange invariant}
l(L) := {\rm dim}\, (V\cap V^{\perp}),
\end{equation}
Note that the intersection pairing is trivial on $V\cap V^{\perp}$.
One immediately gets a useful identity relating these invariants for any $L\subset M$:
\begin{equation} \label{identities}
k(L)+l(L)+s(L) \; =\; {\rm dim}\, (H_n(L)).
\end{equation}

\begin{Definition} \label{def definition P} Let M be a closed oriented $2n$-manifold.
Given a simplicial (or CW) complex $K\subset M$, consider the polynomial
\begin{equation} \label{definition P}
P_{K,M}(X,Y,A,B)=\sum_{L \subset K^{(n)}} {X^{|H_{n-1}(L)|-|H_{n-1}(K)|}Y^{k(L)}A^{s(L)}B^{s^{\perp}(L)}}
\end{equation}
where the sum is taken over all spanning $n$-subcomplexes of $K$.
\end{Definition}

For a detailed discussion of the properties of this polynomial for graphs on surfaces we refer the reader to \cite{Krushkal}. 
We are ready to state the main result of the paper, establishing the duality of the polynomial invariant $P$, generalizing theorem \ref{sphere duality}. Note that the polynomial $T_K$ defined in (\ref{general Tutte}) is a specialization of $P_{K,M}$, see lemma \ref{relation with Tutte} at the end of this section.

\begin{theorem} \label{duality theorem}
{\sl Given a triangulation $K$ of the manifold $M$, let $K^*$ denote the dual cell complex. Then} $$P_{K,M}(X,Y,A,B)=P_{K^*,M}(Y,X,B,A).$$
\end{theorem}

\begin{proof}
The proof consists of two parts, first we establish the duality between the $X$ and $Y$ variables using classical duality theorems from algebraic topology. Then 
duality between the $A$ and $B$ variables will be proved using elements of linear algebra in the presence of a non-degenerate bilinear form (the intersection pairing).
As in the proof of theorem \ref{sphere duality}, for each spanning $n$-subcomplex $L$ of $K$ consider the corresponding ``dual'' spanning $n$-subcomplex $\overline L$ of $K^*$.

Note that for a triangulation $K$ of $M$, $H_{n-1}(K)=H_{n-1}(M)$, therefore the exponent of $X$ in each summand in (\ref{definition P}) equals $|H_{n-1}(L)|-|H_{n-1}(M)|$.

\begin{lemma}$|H_{n-1}(\overline L)|-|H_{n-1}(M)|=k(L)$ \label{identity lemma}
\end{lemma}

\begin{proof}
Consider the homological long exact sequence for the pair $(L,M)$:
\[...\rightarrow H_{n+1}(L) \rightarrow H_{n+1}(M)\rightarrow H_{n+1}(M,L) \rightarrow H_n(L) \rightarrow H_{n}(M) \rightarrow ...\]
Since $L$ does not contain any $(n+1)$-cells, $H_{n+1}(L) = 0$. Recall the following classical theorems of algebraic topology (cf. \cite[Proposition 3.46, Theorem 3.30]{Hatcher}):

\medskip

{\sl  Poincar\'{e}-Lefschetz duality:}  $H_i(M,M\smallsetminus L)\cong H^{n-i}(L)$.
In particular, taking $L=M$, one has

\medskip

{\sl Poincar\'{e} duality:} $H_k(M)\cong H^{n-k}(M)$.

\medskip

Recall from lemma \ref{useful lemma} that $L$ is homotopy equivalent to $M\smallsetminus \overline L$. Then $H_{n+1}(M,L) \cong H^{n-1}(\overline L)$, and $H_{n+1}(M)\cong H^{n-1}(M)$. Also recall that $|H^{n-1}(M)|=|H_{n-1}(M)|$. Coupling these relations with the long exact sequence above gives us that:
\begin{align*}
|H_{n-1}(\overline L)| = {} &|H^{n-1}(\overline L)|  \\
               = {}  &|H_{n+1}(M,L)|\\
               = {} &{\rm rank}\, ({\rm ker}\,(i_*: H_n(L) \rightarrow H_n(M)))+ |H_{n+1}(M)| \\
               = {} & k(L) + |H^{n-1}(M)|\\
               = {} & k(L)+|H_{n-1}(M)|,
\end{align*}

\noindent
concluding the proof of lemma \ref{identity lemma}.
\end{proof}

The following lemma implies that $s(L)=s^{\perp}(\overline L)$, establishing the duality between the $A$ and $B$ variables in the polynomial $P$:

\begin{lemma} \label{v perp lemma}$V(\overline L) \cong V^{\perp}(L)$.
\end{lemma}
\begin{proof}
Decompose $M$ as the union of two submanifolds $\mathcal{H}$ and $\overline{\mathcal{H}}$ that are the handle thickenings of $L$ and $\overline L$ respectively. Denote $\partial:= \partial \mathcal{H}=\partial \overline{\mathcal{H}}$. Suppose $x\in V(\overline L)$. Since the intersection of any $n$-cycle in $\overline{\mathcal{H}}$ with any $n$-cycle in $\mathcal{H}$ is zero, and since $\mathcal{H}$ is a thickening of $L$, it is clear that $x\cdot w=0$ for any $w\in V(L)$. Thus $x\in V^{\perp}(L)$ and so $V(\overline L) \subset V^{\perp}(L)$.

Now we will show that $V^{\perp}(L)\subset V(\overline L).$ Consider the following part of the Mayer-Vietoris sequence:
\begin{equation}\label{MV sequence}\begin{CD}... @>>> H_n(\mathcal{H}) \oplus H_n(\overline{\mathcal{H}}) @>{\alpha}>> H_n(M) @>{\partial}>> H_{n-1}(\partial) @>>> ...
\end{CD}
\end{equation}

Let $x\in H_n(M)$, we claim that $x\notin Im(\alpha)$ implies $x\notin V^{\perp}(L)$. We will establish this by finding an element, $w\in V(L)$ such that $x\cdot w\ne 0$. Since $x\notin Im(\alpha)$ we know by exactness of the above sequence that $x\notin ker(\partial)$ so $\partial (x)=y\in H_{n-1}(\partial)$ is nonzero.
By Poincar\'{e} duality there exists a $z \in H_{n}(\partial)$ such that $y \cdot z\ne 0$. For a moment, we will consider a simpler case (Claim \ref{claim 1}) and then we shall generalize this to the actual problem at hand (Claim \ref{claim 2}).

\begin{claim} \label{claim 1} {\sl Suppose $x_1\in H_n(\mathcal{H},\partial {\mathcal H})$, $y=\partial (x_1)\in H_{n-1}(\partial)$ and $z\in H_n(\partial)$ with $z\cdot y\ne 0$. Then there exists $w\in V(L)$ which pairs nontrivially with $x_1$, i.e. $w\cdot x_1\neq 0$.}
\end{claim}

\begin{proof}
Let $\bar{x}_1 \in H^n(\mathcal{H})$, $\bar{y} \in H^n(\partial)$ and $\bar{z} \in H^{n-1}(\partial)$ be the Poincar\'{e} dual cohomology classes of $x_1$, $y$ and $z$ respectively.  The intersection pairing is the Poincar\'{e} dual of the cup product in cohomology, therefore $\bar{y} \cup \bar{z}\neq 0$.

Our goal is to push the cocycle $\bar{z}$ into $\mathcal{H}$, so that it intersects nontrivially with $x$.  Consider the following commutative diagram in cohomology \cite[8.10]{Dold}:
\[\xymatrixcolsep{5pc}\xymatrix{
H^n(\mathcal{H}) \otimes H^{n-1}(\partial) \ar[r]^{ i^* \otimes Id} \ar[d]^{ Id \otimes \delta} &H^n(\partial) \otimes H^{n-1}(\partial) \ar[r]^{\cup} & H^{2n-1}(\partial) \ar[ld]^{\cong}_{\delta } \\
H^n(\mathcal{H})\otimes H^n(\mathcal{H},\partial)     \ar[r]^{\cup}  &H^{2n}(\mathcal{H},\partial) }
\]
On the level of representatives we have:
\[\xymatrixcolsep{4pc}\xymatrix{
\bar{x}_1 \otimes \bar{z} \ar[r]^{ i^* \otimes Id} \ar[d]^{ Id \otimes \delta} &\bar{y} \otimes \bar{z} \ar[r]^{\cup} &\bar{y} \cup \bar{z} \ar[ld]^{\cong}_{\delta } \\
\bar{x}_1\otimes \delta (\bar{z})     \ar[r]^{\cup}  &\bar{x}_1\cup \delta (\bar{z})}
\]

And since the map on the right is an isomorphism and $\bar{y} \cup \bar{z} \ne 0$, we have that $\bar{x}_1 \cup \delta(\bar{z})\ne 0$. Finally, Poincar\'{e} duality gives an isomorphism $H^n(\mathcal{H}, \partial) \cong H_n(\mathcal{H})$. The image of $\delta(\bar{z})$ under this isomorphism, call it $w$, must pair nontrivially with $x_1$.
\end{proof}

Our situation is slightly different. We have a manifold $M$ decomposed as the union of $\mathcal{H}$ and $\overline{\mathcal{H}}$ along their common boundary. We claim that any cycle $\gamma $ representing a homology class in $H_n(M)$ can be decomposed as the sum of two relative cycles $\gamma _1 + \gamma _2$ in $(\mathcal{H}, \partial)$ and $(\overline{\mathcal{H}}, \partial)$ respectively whose boundaries cancel each other.

If we are able to do this, given $x \in H_n(M)$ satisfying $x\notin Im(\alpha)$, we can decompose $x$ into $x_1+ x_2$ and find a homology class in $H_n(\mathcal{H})$ that intersects either $x_1$ or $x_2$ nontrivially using the method described in Claim \ref{claim 1}.

\begin{claim} \label{claim 2} {\sl  Every cycle $\gamma \in C_n(M)$ is of the form $\gamma = \gamma _1 +\gamma _2$ where $\gamma_1, \gamma_2$ are relative cycles in $(\mathcal{H}, \partial), (\overline{\mathcal{H}}, \partial)$ respectively and $\partial\gamma_1 =-\partial \gamma_2$.}
\end{claim}

\begin{proof}
 Consider the following commutative diagram giving rise to the Mayer Vietoris sequence (\ref{MV sequence}):
\[\xymatrixcolsep{3pc}\xymatrix{
0 \ar[r] &C_n(\partial) \ar[r]^{\phi \phantom{0000}} \ar[d]^{\partial} &{C_n(\mathcal{H}) \oplus C_n(\overline{\mathcal{H}})}  \ar[r]^{\phantom{0000} \psi} \ar[d]^{\partial} & C_n(M) \ar[d]^{\partial } \ar[r] &0 \\
0 \ar[r] &C_{n-1}(\partial) \ar[r]^{\phi \phantom{0000}}  &{C_{n-1}(\mathcal{H}) \oplus C_{n-1}(\overline{\mathcal{H}})} \ar[r]^{\phantom{0000} \psi} & C_{n-1}(M)  \ar[r] &0 }
\]

Let $\gamma $ be a cycle in $C_n(M)$. By the exactness of the top row, there exists $(\gamma _1, \gamma _2) \in C_n(\mathcal{H}) \oplus C_n(\overline{\mathcal{H}})$ such that $\psi(\gamma _1, \gamma _2)= \gamma _1 + \gamma _2 =\gamma $. By definition, $\partial ( \gamma )=0$ since $\gamma$ is a cycle.

Moving to the bottom row, we get that $\psi (\partial \gamma _1,  \partial \gamma _2) = \partial \gamma _1 + \partial \gamma _2= \partial \gamma=0$. Finally we show that $\gamma _1$ and $\gamma _2$ are relative cycles with respect to the pairs $(\mathcal{H}, \partial )$ and $(\overline{\mathcal{H}}, \partial)$. Since $\psi(\partial \gamma _1, \partial \gamma _2)=0$, exactness of the bottom row gives us that there exists some $\gamma_{\partial} \in C_{n-1}(\partial)$ such that $\phi (\gamma_{\partial})= (\partial \gamma _1, \partial \gamma _2)$.
\end{proof}

Now we return to the proof of lemma \ref{v perp lemma}, specifically to the proof of the claim that
$x\notin Im(\alpha)$ implies that $x\notin V^{\perp}(L)$. Suppose $x\notin Im(\alpha)$.
Let $\gamma$ be a cycle representing $x$, then in the notation of claim \ref{claim 2} consider the relative homology class $x_1$ represented by the relative cycle
${\gamma}_1$: $x_1=[{\gamma}_1\in H_n({\mathcal H}, \partial)$. Since $x\notin Im(\alpha)$, $x_1$ and $y=\partial x= \partial x_1$ satisfy the assumptions of claim \ref{claim 1}. Therefore there exists $w\in V(L)$ with $w\cdot x_1=w\cdot x\neq 0$, so $x\notin V^{\perp}(L)$.

Since the goal is to prove $V^{\perp}(L)\subset V(\overline L)$, we may assume
$x \in Im(\alpha)$ and that $x\in V^{\perp}(L)$. Since $Im(\alpha)$ is spanned by $V(L)$, $V(\overline L)$, we may further assume $x\in V(L)$, so $x\in  V(L)\cap V^{\perp}(L)$. Now $x\in V(L)$ means that there is some $x_1\in H_n(\mathcal{H})$ such that $i_{*} (x_1)=x$.
Consider the following part of the long exact sequence:
\[\xymatrix{
.. \ar[r] &H_n(\partial) \ar[r]^{j_*} &H_n(\mathcal{H}) \ar[r]^{k} &H_n(\mathcal{H},\partial) \ar[r] &...
}
\]
We claim that $x_1 \in Im (j_*)$. Suppose to the contrary that $x_1 \notin Im (j_*)$. By exactness this means that $k (x_1)\ne 0$.
Recall that the following intersection pairing is nonsingular \cite{Hatcher}:
\[\xymatrix{
H_n(\mathcal{H}) \times H_n(\mathcal{H},\partial)  \ar[r] & H_{2n}(\mathcal{H}, \partial)
}
\]
So $u\cdot k (x_1) \ne 0$ for some $u\in H_n({\mathcal H})$. Then $u\cdot x=u\cdot x_1= u\cdot k(x_1)\neq 0$, contradicting $x\in V^{\perp}(L)$. So $x_1=j_*(\tilde{x})$ for some $\tilde{x} \in H_n(\partial)$. Recall the Mayer-Vietoris sequence (\ref{MV sequence}):
\[\xymatrix{
& & H_n(\overline{\mathcal{H}}) \ar[dr]^{i_{*}^{'}} & &\\
... \ar[r]& H_n(\partial) \ar[ur]^{j_{*}^{'}} \ar[dr]^{j_*} &\oplus  &H_n(M) \ar[r] &... \\
& & H_n(\mathcal{H}) \ar[ur]^{i_*} & &
}
\]

where ${\alpha}=i_*+i'_*$. Since $x_1 =j_*(\tilde{x})$ and the sequence above is exact, the image $ j_{*}^{'} (\tilde{x})= x_2$ satisfies $i_{*}^{'}(x_2)=x$. Thus $x$ is in $V(\overline L)$, concluding the proof of lemma \ref{v perp lemma}.
\end{proof}
We have shown $V(\overline L)\subset V^{\perp}(L)$ and $V^{\perp}(L) \subset V(\overline L)$, completing the proof of theorem \ref{duality theorem}.
\end{proof}

The proof of theorem \ref{duality theorem} above yields the following statement relating the polynomials $T$ and $P$:

\begin{lemma} \label{relation with Tutte} \sl Let $K$ be a complex embedded in a manifold $M$. The polynomial $T_K$ defined in  (\ref{general Tutte}) is a specialization of $P_{K,M}$:
$$ T_K(X,Y)\; = \; Y^{r/2}\; P_{K,M}(X,Y,Y^{1/2},Y^{-1/2}),$$
where $r={\rm rank}\, (H_n(M))$.
\end{lemma}

Analyzing the Mayer-Vietoris sequence for the decomposition
$M={\mathcal H}\cup \overline{\mathcal{H}}$ establishes the following relation between the invariants
$s, s^{\perp}$ and $l$ of a subcomplex $L\subset M$:
$$s(L)+s^{\perp}(L)+2l(L)\; =\; {\rm rank}\, (H_n(M)).$$
Given this relation and (\ref{identities}), the proof of lemma \ref{relation with Tutte} follows from the fact that the corresponding summands in
the expansions (\ref{general Tutte}), (\ref{definition P}) are equal.

\section{Examples and evaluations} \label{evaluation}

\subsection{Examples.}
In this section we compute the polynomial $P$ for certain handle decompositions of ${\mathbb C}P^2$ and $S^2\times S^2$.
We refer the reader to \cite{Krushkal} for examples of evaluations of $P$ in the case $n=1$ (for graphs on surfaces).

Consider ${\mathbb C}P^2$ with the ``standard'' handle decomposition ${\mathcal H}$, with a single $4$-dimensional $i$-handle $H^i$ for each index
$i=0,2,4$: $${\mathbb C}P^2=H^0\cup H^2\cup H^4.$$
As remarked in section \ref{manifold polynomial}, the polynomial $P_{K,M}$ can be defined not just for a triangulation
$K$ of a manifold $M$ but also in a more general context of a handle decomposition of $M$. (The role of the dual cell complex
$K^*$ is played by the dual handle decomposition.) For the given handle decomposition of ${\mathbb C}P^2$, the sum
(\ref{definition P}) consists of two terms corresponding to $L=H^0$ and $L=H^0\cup H^2$. Since the self-intersection
number of the generating class of $H_2({\mathbb C}P^2)$ is non-trivial, the first term is $B$, and the second term is
$A$. Therefore for this handle decomposition of ${\mathbb C}P^2$ the polynomial is given by
$$P_{{\mathcal H}, {\mathbb C}P^2}=A+B.$$ Observe that the polynomial is the same for a manifold $M$ and for
$M$ with the opposite orientation, in this case $\overline{{\mathbb C}P}^2$, however see the following section for a refinement
that distinguishes them.

Now consider $S^2\times S^2$ with the handle decomposition ${\mathcal H}$ consisting of a single $0$-handle,
two $2$-handles, and a single $4$-handle. The intersection pairing on $H_2(S^2\times S^2)$ is of the form
$\left( \begin{smallmatrix} 0 & 1\\ 1 & 0 \end{smallmatrix} \right)$. There are four summands in the expression
(\ref{definition P}), the term corresponding to $H^0$ (and no $2$-handles) is $B^2$, the term corresponding to $H^0\cup$(a single
$2$-handle) is $1$, the term corresponding to $H^0\cup$(both $2$-handles) is $A^2$, so the polynomial is
$P_{{\mathcal H}, S^2\times S^2}=A^2+2+B^2.$ (In both of these examples the handle decompositions are self-dual and the polynomials
are actually symmetric in $A,B$, giving a stronger version of duality than the general case in theorem \ref{duality theorem}. )

\begin{figure}[ht]
\center{\includegraphics[width=7.5cm]{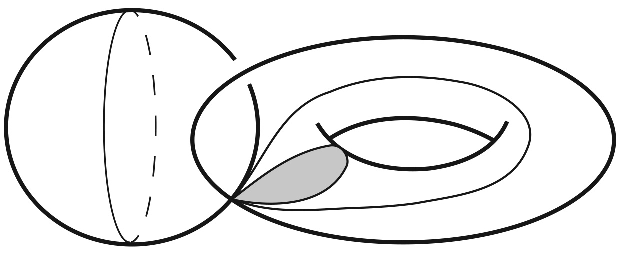}
{\small
    \put(-20,35){$3$}
     \put(-110,23){$1$}
    \put(-203,37){$2$}}}
\caption{A $2$-complex $K$ embedded in ${\mathbb C}P^2$. The labels $1-3$ correspond to the three $2$-cells of $K$.}
\label{K}
\end{figure}

To give an example where all four variables are non-trivially present in the calculation of $P$, consider the $2$-complex $K$ shown in figure \ref{K}. This complex is defined to be the wedge of a $2$-sphere and a $2$-torus, together with a disk bounded by a non-trivial loop on the torus. Therefore $K$ has one $0$-cell, two $1$-cells (generating the first homology of the torus) and three $2$-cells (the top cells of the sphere and the torus, and a $2$-cell attached to a loop on the torus). This complex is embedded in ${\mathbb C}P^2$ so that both the sphere and the torus represent a generator of $H_2({\mathbb C}P^2)$. The terms of the polynomial $K$ corresponding to the subsets of the set 
$\{123\}$ labeling the $2$-cells of $K$ is given in table \ref{K table}. Therefore the polynomial is given by $$P_{K, {\mathbb C}P^2}(X,Y,A,B)\, =\, XB+B+2XA+XYA+2A+YA.$$

\begin{table}[h]
\center{
\begin{tabular}{ 
c | c | c | c | c | c | c | c c}
$\emptyset$  & \{1\} & \{2 \}  & \{3\}  & \{12\}  & \{13\}  & \{23\}  & \{123\}  \\
\hline
$XB$  & $B$  & $XA$ & $XA$ & $A$ & $A$ & $XYA$ & $YA$ & \\
\end{tabular}
\caption{Calculation of the polynomial $P_K$ for the $2$-complex in figure \ref{K}.}
\label{K table}
}
\end{table}

Several variations of the example above may be given, distinguished by the polynomial $P$, but all with the same cellular matroid. For instance, consider the same 
$2$-complex $K$ in figure \ref{K}, but with a different embedding into ${\mathbb C}P^2$. The embedding is defined so that the $2$-sphere represents a generator of $H_2({\mathbb C}P^2)$ and the torus is trivial in second homology of ${\mathbb C}P^2$. In this example the polynomial $P$ equals $$XB+YB+B+XA+XYA+XYB+A+YA.$$
As remarked above and illustrated further by this example, the polynomial $P$ is a strictly stronger invariant than the Tutte polynomial of the cellular matroid, since it depends not just on the complex $K$ but also on its given embedding into a manifold.

The examples in this section are given as just some very basic illustrations of calculations of the polynomial  $P_{K,M}$.

\subsection{Counting simplicial spanning trees} \label{spanning trees}
The value $T_G(1,1)$ of the classical Tutte polynomial equals the number of spanning trees
in a graph $G$, cf. \cite{Bollobas}.
It has been established in \cite{BBC} that the analogous evaluation of the polynomial $T_K$ defined in (\ref{general Tutte}) gives is the number of simplicial spanning trees
of the complex $K$ in the sense of \cite{Kalai}. (The study of simplicial spanning trees has been of considerable recent interest, cf. \cite{DKM, DKM1, Maxwell, MMRW}.)

A {\em weighted} count of spanning trees has been of substantial interest, in particular due to its appearance in the matrix-tree theorem, cf. \cite{DKM}.
Here the weight of an $n$-dimensional spanning tree is the square of the order of its $(n-1)$-st homology group  (which is finite according to the definition of a higher dimensional spanning tree \cite{Kalai}, \cite{DKM}).
Another result of \cite{BBC} is a refinement of the polynomial $T_K(X,Y)$ where the terms in (\ref{general Tutte for triangulation}) are taken with coefficients $|{\rm Tor}(H_{n-1}(L;{\mathbb Z}))|^2$. It is shown in \cite{BBC} that this modified polynomial also satisfies the duality relation
analogous to (\ref{duality equation}), and moreover its evaluation at $(0,0)$ gives the weighted
number of spanning subcomplexes, so it can be calculated by the simplicial matrix-tree theorem. (The evaluation is taken at $(0,0)$ rather than at $(1,1)$ simply due to the fact that the polynomial (\ref{general Tutte}) and its generalization in \cite{BBC} are normalized as the rank-generating polynomial, rather than the Tutte polynomial.)

It is also shown in\cite{BBC} that the Bott polynomial \cite{Bott} of CW complexes may be obtained as a specialization of $T_K(X,Y)$.

\subsection{Other evaluations.} \label{evaluations section}
The polynomial $P_{K, M}(X,Y, A,B)$ reflects both the combinatorial properties of
a complex $K$ and the topological information concerning the embedding of $K$ into $M$.
Recall that the value $T_G(1,1)$ of the classical Tutte polynomial equals
the number of all spanning subgraphs of $G$.
In the following lemma
we point out a generalization of this fact which holds for the polynomial $P$ for graphs on surfaces (corresponding to $n=1$ in definition \ref{def definition P}). Given a graph $G$ embedded in a surface $S$, taking a regular neighborhood of $G$ in $S$ gives it a structure of a {\em ribbon} graph (see \cite{BR, Krushkal} for a detailed account of ribbon graphs). Similarly any subgraph of $G$ then also may be viewed as a ribbon graph.

\begin{lemma} \label{evaluation lemma}  \sl Let $G$ be a graph embedded in an orientable surface $S$. Then
$P_{G,S}(1,1,$ $0,1)$ is the number of (spanning) planar ribbon subgraphs of $G$.
\end{lemma}

The proof is immediate:
$$P_{G,S}(1,1,0,1)=\sum_{H \subset G} 1^{c(H)-c(G)} 1^{k(H)} 0^{s(H)} 1^{s^{\perp}(H)} = |\{ H\subset S:\, s(H)=0.\} |$$
The statement now follows from the fact \cite{Krushkal} that for graphs on surfaces $s(H)$ is twice the genus
of the regular neighborhood of the graph $H$ in $S$.
\qed

Spanning {\em quasi-trees} have recently been used by several authors as analogues of spanning trees that are suitable in the context of ribbon graphs, cf. \cite{CKS, Butler}. Recall that a quasi-tree is defined as a connected subgraph $H$ of a ribbon graph $G$ such that the boundary of $H$ is connected. Reformulating this definition in terms of the invariants introduced at the beginning of section \ref{manifold polynomial} and generalizing to higher dimensions, define a {\em simplicial spanning quasi-tree} of $K\subset M$ to be a
spanning $n$-subcomplex $L$ such that ${\rm dim} H_{n-1}(L)=0$, and $ k(L)=l(L)=0$. For graphs this definition coincides with the one discussed above.
To ensure that a given complex $K$ embedded in $M$ has a simplicial spanning quasi-tree we require that the $(n-1)$-st Betti number ${\beta}_{n-1}(K)=0$, this is related to the notion of a complex which is acyclic in positive codimension \cite{DKM, BBC}. (For example, one requires that a graph is connected to ensure that a spanning tree exists.)

A slight variation of the definition (\ref{definition P}) is the following $5$-variable polynomial:
\begin{equation} \label{definition P'}
\widetilde P_{K,M}(X,Y,A,B,C)=\sum_{L \subset K^{(n)}} {X^{|H_{n-1}(L)|-|H_{n-1}(K)|}Y^{k(L)}A^{s(L)}B^{s^{\perp}(L)}C^{l(L)}}
\end{equation}

Using the relation (\ref{identities}) between the parameters $k, l, s$ it is easy to see that the polynomials $P, \widetilde P$ carry equivalent information.

\begin{lemma} \label{quasitree lemma}  \sl $\widetilde P_{K,M}(0,0,1,1,0)$ is the number of spanning quasi-trees of a complex $K$ embedded in a manifold $M$.
\end{lemma}

It seems likely that the results of \cite{Butler} can be extended from graphs to higher-dimensional simplicial complexes to obtain a simplicial quasi-tree expansion of the polynomial $P$.

\section{Generalizations of the polynomial invariants ${T}$, $P$.} \label{generalizations}

\subsection{Polynomials ${\mathbf{{T}^j}}$ for triangulations of the sphere ${\mathbf S^N}$.} \label{other T}

In this subsection we note that definition \ref{sphere definition} of the polynomial $T$ may be extended to
spheres of any (not necessarily even) dimension, giving rise to a collection of polynomials ${T}^j$:

\begin{Definition}
Let $K$ be a simplicial complex of dimension $n$ and let $1\leq j \leq n$. Consider the polynomial invariant
\[ T^j_K(X,Y)=\sum_{L \subset K^{(j)}} {X^{|H_{j-1}(L)|-|H_{j-1}(K))}Y^{|H_j(L)|}}\]
where $K^{(j)}$ is the $j$-skeleton of $K$ and the summation is taken over all spanning $j$-subcomplexes $L$ of $K$ such that
$L^{(j-1)}=K^{(j-1)}$.
\end{Definition}

Clearly the original polynomial $T$ in definition \ref{sphere definition} (for $N=2n$) equals ${T}^n$ in the definition above.
The analogue of theorem \ref{sphere duality} for the polynomials ${T}^j$ is stated as follows:

\begin{lemma} \label{sphere duality k}
{\sl Given a triangulation $K$ of $S^{N}$, let $K^*$ denote the dual cell complex. Then
$$T^j_K(X,Y)=T^{N-j}_{K^*}(Y,X)$$
for $1\leq j\leq N$.}
\end{lemma}

\subsection{The polynomial $\overline P $ for triangulations of oriented $4n$ dimensional manifolds.}

When the dimension of an oriented manifold $M$ is divisible by $4$, the definition of the polynomial $P$ may be refined
further. Following the notation used in equations (\ref{s}), observe that the intersection pairing is a symmetric non-degenerate
bilinear form on the vector spaces
$V/(V\cap V^{\perp})$, $V^{\perp}/(V\cap V^{\perp})$. Denote by $s_+(L)$ the dimension of a maximal subspace of
$V/(V\cap V^{\perp})$ on which the intersection pairing is positive definite, and similarly by $s_-(L)$ the dimension
where it is negative definite. $s^{\perp}_+(L)$, $s^{\perp}_+(L)$ are defined analogously. Note that
$$s(L)=s_+(L)+s_-(L), \; \, s^{\perp}(L)=s^{\perp}_+(L)+s^{\perp}_-(L).$$
Given a triangulation $K$ of $M^{2n}$, where $n$ is even, consider

\begin{equation} \label{definition overline P}
\overline P_{K,M}(X,Y,A_+, A_-,B_+, B_-)\; =
\end{equation}
$$=\; \sum_{L \subset K^{(n)}}
{X^{|H_{n-1}(L)|-|H_{n-1}(M)|}Y^{k(L)}A_+^{s_+(L)}A_-^{s_-(L)} B_+^{s_+^{\perp}(L)}B_-^{s_-^{\perp}(L)}}
$$

This is a refinement of the polynomial (\ref{definition P}) in the sense that
$$P_{K,M}(X,Y,A,B) = \overline P_{K,M}(X,Y,A, A,B, B).$$
Note that while reversing the orientation of the manifold $M$ did not change the polynomial
$P$, the polynomial $\overline P$ changes as follows:
$$\overline P_{K,\overline M}(X,Y,A_+, A_-,B_+, B_-) = \overline P_{K,M}(X,Y,A_-, A_+,B_-, B_+),$$
where $\overline M$ denotes $M$ with the opposite orientation.
The duality theorem \ref{duality theorem} takes the form
$$\overline P_{K,M}(X,Y,A_+, A_-,B_+, B_-)=\overline P_{K^*,M}(Y,X,B_+, B_-,A_+, A_-).$$

\section{Remarks and questions} \label{questions}

We conclude by listing several questions motivated by our results.

{\bf 1.} An example of an evaluation of the polynomial $P$ is given in section \ref{evaluations section}: $P_{G,S}(1,1,$ $0,1)$ is the number
of planar ribbon subgraphs of a graph $G$ embedded in a surface $S$. The results of \cite{BBC}, discussed in section \ref{spanning trees},
show that the evaluation of $T_K(0,0)$ is the number of  simplicial spanning trees of a complex $K$.
It is likely that there are other evaluations of these polynomials which reflect both the
combinatorics of the triangulation and the topology of the ambient manifold $M$. For example, an interesting question is whether a higher-dimensional generalization
of lemma \ref{evaluation lemma}  holds: Given a complex $K\subset M^{2n}$, is $P_{K,M}(1,1,0,1)$ equal to the number of those subcomplexes of $K$ whose
neighborhoods in $M$ embed in a homology $2n$-dimensional sphere?

{\bf 2.} The chromatic polynomial (which may be thought of as a one-variable specialization of the Tutte polynomial) of planar graphs
is known to satisfy a sequence of linear local relations when it is evaluated at the Beraha numbers, and moreover
it satisfies a remarkable quadratic {\em golden identity} at the golden ratio \cite{Tutte3, FK}. It is a natural generalization of this problem to ask whether there are evaluations of the polynomials introduced in this paper which satisfy additional local
relations. The planar identities satisfied by the Tutte polynomial are known to fit in the
framework of quantum topology (see \cite{FK}), and the lack of interesting topological quantum field theories
in higher dimensions indicates that possible local relations in higher dimensions would have to be of a different nature.

{\bf 3.} A generalization of the Tutte and of the Bollob\'{a}s-Riordan polynomials, motivated by ideas in quantum gravity, in the context of tensor graphs has been introduced in  \cite{Tanasa}. It would be interesting to find out if there is a relation between the invariants defined in this paper and those in \cite{Tanasa}.

{\bf 4.} Recall that the polynomial $T_K$ defined in (\ref{general Tutte}) may be formulated in the context of matroid theory (section \ref{matroid section}). Such a formulation of the more general polynomial $P_{K,M}$ defined in (\ref{definition P}) is not immediate. A different relation between the polynomial P (in the context of graphs on surfaces) and matroids is presented in \cite{ACEMS}. It seems reasonable that this approach (due to Las Vergnas \cite{Las Vergnas}, in terms of {\em matroid perspectives}) may generalize to higher dimensions as well.


\medskip


\begin{thebibliography}{10}

\bibitem{ACEMS} R. Askanazi, S. Chmutov, C. Estill, J. Michel and P. Stollenwerk,
{\em Polynomial invariants  of graphs on surfaces},   Quantum Topol. 4 (2013), 77-90. [arXiv:1012.5053]

\bibitem{BBC} C. Bajo,  B. Burdick, S. Chmutov, {\em On the Tutte-Krushkal-Renardy polynomial for cell complexes}, J. Combin. Theory Ser. A 123 (2014), 186-201.
[arXiv:1204.3563]

\bibitem{Bollobas} B. Bollob\'{a}s,  Modern graph theory, Springer, 1998.

\bibitem{BR} B. Bollob\'{a}s and O. Riordan, {\em A polynomial invariant of graphs on orientable surfaces}, Proc. London Math. Soc. 83(2001), 513-531.

\bibitem{Bott} R. Bott, {\em Two new combinatorial invariants for polyhedra}, Portugaliae Math. 11 (1952), 35-40. 

\bibitem{BO}  T. Brylawski and J.  Oxley, {\em The Tutte polynomial and its applications},  Matroid applications, 123-225, Encyclopedia Math. Appl., 40, Cambridge Univ. Press, Cambridge, 1992. 

\bibitem{Butler} C. Butler, {\em A quasi-tree expansion of the Krushkal polynomial},
arXiv:1205.0298.

\bibitem{CKS} A. Champanerkar, I. Kofman and N. Stoltzfus, {\em Quasi-tree expansion for the Bollob\'{a}s-Riordan-Tutte polynomial}, Bull. Lond. Math. Soc. 43 (2011), 972-984.

\bibitem{Cordovil} R. Cordovil\ and B. Lindstr\"{o}m, {\em Simplicial matroids}, Combinatorial geometries, 98-113, Encyclopedia Math. Appl., 29, Cambridge Univ. Press, Cambridge, 1987

\bibitem{CrapoRota} H. Crapo and G.-C. Rota,
{\em Simplicial geometries}, Combinatorics (Proc. Sympos. Pure Math., Vol. XIX), 71-75. Amer. Math. Soc., Providence, R.I., 1971

\bibitem{Dold} A. Dold,
Lectures on Algebraic Topology, Berlin: Springer, 1972.

\bibitem{DKM} A.M. Duval, C.J. Klivans and J. Martin, {\em Simplicial matrix-tree theorems},
 Trans. Amer. Math. Soc. {\bf 361}(11) (2009) 6073--6114.

\bibitem{DKM1} A.M. Duval, C.J. Klivans and J.L. Martin, {\em Cellular spanning trees and Laplacians of cubical complexes}, Adv. in Appl. Math. 46 (2011), 247-274.

\bibitem{FK} P. Fendley and V. Krushkal, {\em Tutte chromatic identities from the Temperley-Lieb algebra},
Geom. Topol. 13 (2009), 709-741 [arXiv:0806.3484]

\bibitem{Hatcher} A. Hatcher, {Algebraic Topology}, Cambridge: Cambridge University Press, 2009.

\bibitem{Kalai} G. Kalai, {\em Enumeration of ${\mathbb Q} $-acyclic simplicial complexes},
Isr. J. Math. 45 (1983), 337-351.

\bibitem{Krushkal} V. Krushkal, {\em Graphs, links, and duality on surfaces}, Combin. Probab. Comput. 20 (2011), 267-287
[arXiv:0903.5312]


\bibitem{Las Vergnas} M. Las Vergnas, On the Tutte polynomial of a morphism of matroids, Ann. Discrete Math. 8
(1980), 7-20.


\bibitem{MMRW} J.L. Martin, M. Maxwell, V. Reiner and Scott O. Wilson,
{\em Pseudodeterminants and perfect square spanning tree counts}, 
arXiv:1311.6686.

\bibitem{Maxwell} M. Maxwell, {\em  Enumerating bases of self-dual matroids}, J. Combin. Theory Ser. A 116 (2009), 351-378.

\bibitem{Milnor} J. Milnor,  Morse theory, Princeton University Press, 1963.

\bibitem{Oxley}  J.G. Oxley, Matroid theory, Oxford University Press, New York, 1992.

\bibitem{RS} C. Rourke and B. Sanderson, {Introduction to Piecewise-linear Topology}, Berlin: Springer-Verlag, 1972.


\bibitem{Tanasa} A. Tanasa, {\em Generalization of the Bollob\'{a}s-Riordan polynomial for tensor graphs}, J. Math. Phys. 52, 073514 (2011).

\bibitem{Tutte1} W.T. Tutte, {\em A contribution to the theory of chromatic polynomials}, Canadian J. Math.
6 (1954), 80-91.


\bibitem{Tutte3} W.T.~Tutte, {\em On chromatic polynomials and the golden ratio}, J. Combinatorial Theory 9 (1970),
289-296.


\bibitem{Tutte4} W.T. Tutte, Introduction to the theory of matroids,
American Elsevier Publishing Co., Inc., New York 1971

\bibitem{Welsh} D. J. A. Welsh,  Matroid theory,  Academic Press, London-New York, 1976.

\bibitem{Whitney} H. Whitney, {\em On the abstract properties of linear dependence}, Amer. J. Math. 57 (1935), 509-533.

\end{thebibliography}
\end{document}